\documentclass[11pt]{amsart}

\usepackage[square,compress,comma, numbers,sort]{natbib}
\usepackage[colorlinks=true, citecolor=red, linkcolor=blue]{hyperref}
\usepackage{amsfonts,mathtools}

\allowdisplaybreaks[4]

\usepackage{amsmath}
\usepackage{bbm}
\usepackage{amssymb}
\usepackage{mathtools}

\usepackage{color}

\def\zE#1{{\textcolor{c30}{#1}}}

\definecolor{c20}{rgb}{0.,0.7,0.}
\definecolor{c30}{rgb}{0.,0.,1.}
\definecolor{c40}{rgb}{1,0.1,0.7}
\definecolor{c50}{rgb}{1,0,0}
\definecolor{c60}{rgb}{1,0.9,0.1}

\newcommand{\E}[1]{\mathbb{E}\left\{ #1\right\}}

\newcommand{\pk}[1]{\mathbb{P} \left\{ #1 \right \} }

\newcommand{\R}{\mathbb{R}}

\newcommand{\N}{\mathbb{N}}

\newcommand{\limit}[1]{\lim_{#1 \to \infty}}

\newcommand{\BQN}{\begin{eqnarray}}
\newcommand{\EQN}{\end{eqnarray}}
\newcommand{\BQNY}{\begin{eqnarray*}}
\newcommand{\EQNY}{\end{eqnarray*}}

\newcommand{\BS}{\begin{sat}}
\newcommand{\ES}{\end{sat}}
\newcommand{\BT}{\begin{theo}}
\newcommand{\ET}{\end{theo}}
\newcommand{\BK}{\begin{korr}}
\newcommand{\EK}{\end{korr}}

\newcommand{\BD}{\begin{de}}
\newcommand{\ED}{\end{de}}

\newcommand{\BIT}{\begin{itemize}}
\newcommand{\EIT}{\end{itemize}}
\newcommand{\BDI}{\begin{description}}
\newcommand{\EDI}{\end{description}}

\newcommand{\BRM}{\begin{remarks}}
\newcommand{\ERM}{\end{remarks}}

\newcommand{\BEL}{\begin{lem}}
\newcommand{\EEL}{\end{lem}}

\newtheorem{theo}{Theorem}[section]
\newtheorem{sat}[theo]{Proposition}
\newtheorem{de}[theo]{Definition}
\newtheorem{lem}[theo]{Lemma}

\newtheorem{korr}[theo]{Corollary}
\newtheorem{remark}[theo]{Remark}
\newtheorem{remarks}[theo]{Remarks}
\newtheorem{prop}[theo]{Proposition}

\newcommand{\nelem}[1]{{Lemma \ref{#1}}}

\newcommand{\netheo}[1]{{Theorem \ref{#1}}}

\newcommand{\prooftheo}[1]{ \textsc{\bf Proof of Theorem} \ref{#1}:}
\newcommand{\proofprop}[1]{\textsc{\bf Proof of Proposition} \ref{#1}:}
\newcommand{\prooflem}[1]{\textsc{\bf Proof of Lemma} \ref{#1}:}

\newcommand{\COM}[1]{}

\newcommand{\auxtheo}[1]{ \textsc{\bf Auxiliary results for Theorem} \ref{#1}:}

\def\td{\text{\rm d}}
\def\IF{\infty}

\newcommand{\QED}{\hfill $\Box$}

\def\IF{\infty}

\def\bqny#1{{\begin{eqnarray*} #1 \end{eqnarray*}}}
\def\bqn#1{{\begin{eqnarray} #1 \end{eqnarray}}}

\begin{document}

\title{Parisian \& Cumulative Parisian Ruin Probability for Two-Dimensional Brownian Risk Model}

\author{Nikolai Kriukov }
\address{Nikolai Kriukov, Department of Actuarial Science, 
University of Lausanne,\\
UNIL-Dorigny, 1015 Lausanne, Switzerland
}
\email{nikolai.kriukov@unil.ch}

\bigskip

\date{\today}
\maketitle

{\bf Abstract:} Parisian ruin probability in the classical Brownian risk model, unlike the standard ruin probability can not be explicitly calculated even in one-dimensional setup. Resorting on asymptotic theory, we derive in this contribution an asymptotic approximations of both Parisian and cumulative Parisian ruin probability and simultaneous ruin time for the two-dimensional Brownian risk model when the initial capital increases to infinity.



\section{Introduction} 
Calculation of Parisian ruin for Brownian risk model has been initially considered in \cite{LCPalowski}. For general Gaussian risk models Parisian ruin cannot be calculated explicitly. As shown in \cite{MR3457055,MR3414985} methods from the theory of extremes of Gaussian random fields can be successfully applied to approximate the Parisian ruin for general Gaussian risk models. In this paper, we shall be concerned with the classical bivariate Brownian motion risk model, which in view of 
recent findings in \cite{delsing2018asymptotics}, appears naturally as the limiting model of some general bivariate insurance risk model. Consider therefore two insurance risk portfolios with corresponding risk models 
\[R_1(t)= u+ c_1t - W_1(t), \quad R_2(t)= au+ c_2t - W_2(t), \quad t\ge 0,\]
where $W_1,W_2$ are two standard Brownian motion and the initial capital for the first portfolio is $u>0$, whereas for the second it is equal $au, a\in (-\IF, 1]$. Further $c_1$ and $c_2$ denote 
the premium rates of the first and the second portfolio, respectively. In this contribution we shall consider the benchmark model where 
$(W_1(t),W_2(t)),t\ge 0$ are assumed to be jointly Gaussian with the same law as
\BQN \label{BB}
(B_1(t), \rho B_1(t)+ \sqrt{1- \rho^2} B_2(t)), \quad t\ge 0, \quad \rho \in (-1,1),
\EQN
where $B_1,B_2$ are two independent standard Brownian motions. As mentioned above, this model is supported by the findings of \cite{delsing2018asymptotics}.

For given $T>0$ and $H\ge 0$ define the simultaneous Parisian ruin probability on finite time horizon $[0,T]$
\bqn{ 
p_{H,T}(u,au) = \pk{ \exists t\in [0,T], \forall s \in [t,t+H]:\ R_1(s) < 0, R_2(s) < 0}.
} 
When $H=0$, the simultaneous Parisian ruin reduces to simultaneous classical ruin studied recently in \cite{mi:18}

It follows easily that the Parisian ruin probability is smaller than the simultaneous ruin probability $p_{0,T}(u,au)$, namely 
$$ p_{H,T}(u,au) \le p_{0,T} = \pk{ \exists t\in [0,T]: R_1(t)< 0, R_2(t) < 0 }$$
for any $u,H,T$ positive since $p_{H,T}(u,au)$ is monotone in $H$. In \cite{mi:18} it is shown that 
the simultaneous ruin probability can be simply bounded as follows 
\begin{eqnarray}
\label{bXh}
\pk{W_1^*(T)>u ,W_2^*(T)>au}
\le
p_{0,T}(u,au)\le
\frac{\pk{W_1^*(T)>u,W_2^*(T)>au }}{\pk{W_1(T)>\max(c_1 T,0),W_2(T)>\max(c_2T,0)}},
\end{eqnarray}
where we set $W_i^*(t)=W_1(t)- c_i t, i=1,2.$ \\
A simple lower bound for $p_{H,T}(u, au)$ is the following 

\bqn{ \pk{ \forall t\in [T,T+H]: R_1(t)< 0, R_2(t) < 0 } \le p_{H,T}(u,au), \label{4}
}
which is valid for any $u>0$. The above lower bound is very difficult to evaluate even asymptotically when $u$ tends to infinity. A simpler case is when $a < \rho, \rho >0$. We have (see Appendix) that for all large $u$ and some 
$C \in (0,1)$ 
\bqn{ \quad C\pk{ \forall t\in [T,T+H]: W_1^*(t)> u } \le p_{H,T}(u,au) \le \pk{ \sup_{ t\in [0,T]}W_1^*(t)> u }.\label{5}
}
Since $\pk{ \sup_{ t\in [0,T]}W_1^*(t)> u }$ can be evaluated explicitly, it follows easily that it is asymptotically as $u\to \IF$ equal to 
$\pk{ W_1^*(T)> u}$ and by \cite{MR3457055}[Thm 2.1] the lower bound is proportional to $\pk{ W_1^*(T)> u}/u$ as $u\to \IF$. 
Therefore, even for this simple case, the bounds derived above do not capture the exact decrease of the Parisian ruin probability as $u\to \IF$. The reason for this is that the interval $[T,T+H]$ is quite large. In the sequel, under the restriction that $H= S /u^2$ for any $S\ge 0$ we show that it is possible to derive the exact approximations of the Parisian ruin probability.

Motivated by \cite{KrzyszSjT} in this paper we shall also investigate the so called cumulative Parisian ruin probability on the finite time interval $[0,T]$, i.e., 
$$ \Psi_{L,T}(u,au)=\pk{ \int_{0}^T \mathbb{I}( R_1(t)< 0, R_2(t) < 0 ) dt > L/f(u)},$$
where $L>0$ is a given constant and $f(u)$ is some positive function depend on $u$. It is clear that the above is bounded by 
$p_{0,T} (u, au)$ and the calculation of the cumulative Parisian ruin probability is not possible for any fixed $u$ and $x$ positive. A natural question here is (see \cite{KrzyszSjT} for the infinite time-horizon case) if we can approximate the cumulative Parisian ruin probability when $u$ tends to infinity. This in particular requires to determine explicitly the function $f(u)$. In the case of one-dimensional risk model it is shown recently in \cite{KrzyszSjT} that the cumulative Parisian ruin probability (or in the language of that paper the tail of the sojourn time/occupation time) can be approximated exactly when $u$ becomes large. In that aforementioned paper $f(u)$ equals $u^2$. We shall show that this is the right choice also for our setup. 

Brief organization of the paper. Next section presents the exact asymptotics of both Parisian and cumulative Parisian ruin as $u\to \IF$. Additionally, we discuss the approximation of the cumulative Parisian ruin time $\tau_L(u)=\inf_{T>0}\int_0^{T}\mathbb{I}(R_1(t)<0,R_2(t)<0)\td t>L/f(u)$. 

Section 3 is dedicated to the proofs. We conclude this contribution with an Appendix.

\section{Main Results}

For simplicity we consider below the case $T=1$ and $a\le 1$, since by the self-similarity of Brownian motion we can resolve the other portfolio. Let in the following
\BQN \label{lam12}
\lambda_1= \frac{1- a\rho}{1- \rho^2}, \quad \lambda_2=\frac{a- \rho }{1- \rho^2},
\EQN
which are both positive if $a \in (\rho, 1]$. For the particular choice of $H= S/u^2$ we shall denote $P_{H}(u,au)$ simply as 
$\psi_{S}(u, au)$. We consider first the approximation of the Parisian ruin.

\BT \label{Th1} Let $c_1,c_2$ be two given real constants and $S\ge 0$ be given. \\
i) If $a\in (\rho,1]$, then as $u\to \IF$
\BQN
\psi_{S}(u, au) \sim C_{a,\rho}(S)\pk{W_1^*(1)>u,W_2^*(1)>au},
\EQN
where the constant $C_{a,\rho}(S)\in (0,\IF)$ is given by 
\BQN
\quad C_{a,\rho}(S)= \lambda_1\lambda_2\int_{\R^2}
\pk{\exists t \ge 0 , \forall s \in [t-S,t]:
\begin{array}{ccc}
W_1(s) - s > x \\
W_2(s) - as > y
\end{array}
} e^{ \lambda_1 x+ \lambda_2 y }\, \td x \td y.
\label{car}
\EQN

ii) If $a\le \rho$, then we have as $u\to \IF$
\BQN
\psi_{S}(u, au) \sim C(S)\pk{W_1^*(1)>u,W_2^*(1)>au},
\EQN
where
$C(S)= \E{ e^{ \sup_{t \ge 0} \inf_{s \in [t-S,t]} ( W_1(s)- s)}}< \IF.$ 
\ET

The approximation of the cumulative Parisian ruin requires some different arguments since the sojourn functional is different from the supremum functional. In the following we shall choose the scaling function $f(u)$ to be equal to $u^2$. Since $T=1$, we can omit it and write simply $\Psi_{L}(u,au)$ instead of $\Psi_{T,L}(u,au)$.

\BT \label{Th2} Under the setup and the notation of \netheo{Th1} for any $L>0$ we have:\\
i) If $a\in (\rho,1]$, then as $u\to \IF$
\BQN
\Psi_{L}(u, au) \sim K_{a,\rho}(L) \pk{W_1^*(1)>u,W_2^*(1)>au},
\EQN
where the constant $K_{a,\rho}(L)\in (0,\IF)$ is given by 
\BQN
\quad K_{a,\rho}(L)=\lambda_1\lambda_2\int_{\R^2}\pk{\int_{0}^{\infty}\mathbb{I}(W_1(t)-t>x,W_2(t)-at>y)\td t>L}e^{\lambda_1 x+\lambda_2 y}\td x\td y.
\label{car2}
\EQN

ii) If $a\le \rho$, then we have as $u\to \IF$
\BQN
\Psi_{L}(u, au) \sim K(L)\pk{W_1^*(1)>u,W_2^*(1)>au},
\EQN
where 
\bqn{
K(L)&=&\int_{\R}e^x\pk{\int_0^{\infty}\mathbb{I}(W_1(t)-t>x)\td t>L}\td x \in(0,\infty). \label{Ka}
}
\ET
\begin{remark} Theorems \ref{Th1} and \ref{Th2} may be used also if $c_1$, $c_2$, $S$ and $L$ depend on $u$, but have finite limits as $u\to\infty$ ($c_1(u)\to c_1^*$, $c_2(u)\to c_2^*$, $S(u)\to S^*$, $L(u)\to L^*$). In this case all constants $S$ and $L$ on the right-hand sides should be replaced by $S^*$ and $L^*$, respectively.
\end{remark}
The asymptotic distribution of cumulative ruin time $\tau_L(u)$ may be explicitly calculated from Theorem \ref{Th2} by using the self-similarity of Brownian motion.

\begin{prop}\label{Th3}
i) If $a\in(\rho,1]$, then for any $x\in(0,\infty)$ and $0\leq L_2\leq L_1\leq 1$ 
\bqny{
\lim_{u\to\infty}\pk{u^2(1-\tau_{L_1}(u))\geq x|\tau_{L_2}(u)\leq 1}= \Gamma(L_1,L_2)e^{-x\frac{1-2a\rho+a^2}{2-2\rho^2}},
}
where with $K_{a,\rho}$ defined in $\eqref{car2}$ we have
\bqny{
\Gamma(L_1,L_2)=\frac{K_{a,\rho}(L_1)}{K_{a,\rho}(L_2)}.
}
ii) If $a\leq\rho$, then for any $x\in(0,\infty)$ and $0\leq L_2\leq L_1\leq 1$
\bqny{
\lim_{u\to\infty}\pk{u^2(1-\tau_{L_1}(u))\geq x|\tau_{L_2}(u)\leq 1}=\Gamma(L_1,L_2)e^{-\frac{x}{2}},
}
where with $K$ defined in $\eqref{Ka}$ we have
\bqny{
\Gamma(L_1,L_2)=\frac{K(L_1)}{K(L_2)}.
}
\end{prop}

\section{Proofs}

\auxtheo{Th1} Let in the following $\delta(u,T)=1-Tu^{-2}$ for $T,u>0$.

Note that for any $S,T$ positive 
\bqn{
m(u,S,T)&\coloneqq & \pk{\exists_{t \in [0,\delta(u,T)]}, \forall s \in [0, S/u^2]: W_1^*(s)> u, W_2^*(s)> au } \notag
\\
& \le & \pk{\exists_{t \in [0,\delta(u,T)]}: W_1^*(t)> u, W_2^*(t)> au } \notag
\\
&\le & e^{-T/{8}}
\frac{\pk{W_1^*(1)\ge u, W_2^*(1)\ge au}}{\pk{W_1(1)>\max(c_1,0),W_2(1)>\max(c_2,0)}}
}
for all $u$ large, where the upper bound follows from \cite{mi:18}[Lemma 4.1] (see Appendix, Lemma $\ref{inimp}$).

We give below the exact asymptotics of 
\[
M(u,S,T)\coloneqq\pk{\exists_{t \in [\delta(u,T),1]}, \forall s \in [t, t+S/u^2]: W_1^*(s)> u, W_2^*(s)> au }
\]
as $u$ tends to infinity.

\BEL
i) For any $a\in (\rho, 1]$ and any $S,T>0$ we have as $u\to\infty$
\label{singA}
\BQN
M(u,S, T)&\sim & u^{-2} \varphi_\rho(u+c_1,au+c_2) I(S,T) ,
\EQN
where
$$ I(S,T)\coloneqq \int_{\R^2}
\pk{t\in [0, T ], \forall s\in [t-S,t]
\begin{array}{ccc}
W_1(s) - s > x \\
W_2(s) - as > y
\end{array}
} e^{ \lambda_1 x + \lambda_2 y }\, \td x \td y \in (0,\IF).$$
ii) For any $a \le \rho$ and any $S,T>0$ we have as $u\to\infty$
\BQNY
M(u,S,T) \sim u^{-1} \varphi_\rho(u+ c_1,\rho u+ c_2) I(S,T),
\EQNY
where
$$I(S,T)\coloneqq\int_{\R^2}\pk{ \exists t\in [0,T], \forall s\in [t-S,t]: W_1(s)-s> x} \Bigl[\mathbb{I}(a <\rho)+ \mathbb{I}(y\zE{<} 0,a=\rho)\Bigr] e^{ x - \frac{y^2-2{y(c_2- c_1 \rho)}}{2 (1- \rho^2)}}\, \td x \td y. $$
\label{L1}
\EEL

\prooflem{L1}
i) For any $x,y\in \R$ put
\bqny{
u_x=u+c_1-x/u,\qquad u_y=au+c_2-y/u.
}
Writing $\varphi(x,y)$ for the joint pdf of $(W_1(1),W_2(1))^\top$ we have
\bqn{
\varphi_{\rho}(u_x,u_y)=:\varphi_\rho(u+c_1,au+c_2)\psi_{u}(x,y),\label{psi_u}
}
where (write $\Sigma$ for the covariance matrix of $(W_1(1),W_2(1))^\top$)
\bqn{
\begin{aligned}
\log\psi_u(x,y)&=\frac{1}{u^2}(u+c_1,au+c_2)\Sigma^{-1}(x,y)^\top-\frac{1}{2u^2}(x,y)\Sigma^{-1}(x,y)^\top\\
&=(1,a)\Sigma^{-1}(x,y)^\top+\frac{1}{u}(c_1,c_2)\Sigma^{-1}(x,y)^\top-\frac{1}{2u^2}(x,y)\Sigma^{-1}(x,y)^\top\\
&\to(1,a)\Sigma^{-1}(x,y)^\top,\qquad\qquad\qquad\qquad\qquad\qquad u\to\infty\\
&=\frac{1}{1-\rho^2}(1,a)\begin{pmatrix}1 & -\rho\\-\rho & 1\end{pmatrix}(x,y)^\top\\
&=\frac{1-a\rho}{1-\rho^2}x+\frac{a-\rho}{1-\rho^2}y=\lambda_1x+\lambda_2y.
\end{aligned}\label{ipsi}
}
Denote further
\bqny{
u_{x,y}=u_y-\rho u_x=(a-\rho)u-(y-\rho x)/u+ c_2-\rho c_1.
}
Let $B_1,~B_2$ be two independent Brownian motions. The representation of $(W_1(t),W_2(t))$ in terms of $B_1$ and $B_2$ is given in $\eqref{BB}$.
Define the following time transform
\bqny{
\bar{s}_u=1-s/u^2,
}
and set
\bqny{
A(u)=u^{-2}\varphi_{\rho}(u+c_1,au+c_2).
}
For the function $M(u,S,T)$ we have with $\psi_u$ defined in \eqref{psi_u}
\bqny{
& &M(u,S,T)\\
& &\quad=u^{-2}\int_{\R^2}\pk{\exists t\in[\delta(u,T),1]\forall s\in[t,t+S/u^2]:\begin{split}&W_1^*(s)>u\\&W_2^*(s)>au\end{split}\left|\begin{split}&W_1(1)=u_x\\&W_2(1)=u_y\end{split}\right.}\varphi_\rho(u_x,u_y)\td x\td y\\
& &\quad=A(u)\int_{\R^2}\pk{\exists t\in[0,T]\forall s\in[t-S,t]:\begin{split}&W_1^*(\bar{s}_u)>u\\&W_2^*(\bar{s}_u)>au\end{split}\left|\begin{split}&W_1(1)=u_x\\&W_2(1)=u_y\end{split}\right.}\psi_u(x,y)\td x\td y\\
& &\quad=A(u)\int_{\R^2}\pk{\exists t\in[0,T]\forall s\in[t-S,t]:\begin{split}&B_1(\bar{s}_u)-c_1\bar{s}_u>u\\&\rho B_1(\bar{s}_u)+\rho^* B_2(\bar{s}_u)-c_2\bar{s}_u>au\end{split}\left|\begin{split}&B_1(1)=u_x\\&\rho^*B_2(1)=u_{x,y}\end{split}\right.}\psi_u(x,y)\td x\td y\\
& &\quad=:A(u)\int_{R^2}h_u(T,S,x,y)\psi_u(x,y)\td x\td y.
}
Define two auxiliary processes for $s\in[-S,T]$ as follows
\bqn{
\begin{split}
B_{u,1}(s)&:=\left\{B_1(\bar{s}_u)\left|B_1(1)=u_x\right.\right\}-\bar{s}_uu_x,\\
B_{u,2}(s)&:=\left\{B_2(\bar{s}_u)\left|\rho^* B_2(1)=u_{x,y}\right.\right\}-\bar{s}_uu_{x,y}/\rho^*.
\end{split}\label{auxBB}}
We can represent the function $h_u(T,S,x,y)$ in terms of these processes as follows
\bqny{
h_u(T,S,x,y)=\pk{\exists t\in[0,T]\forall s\in[t-S,t]:\begin{split}&u(B_{u,1}(s)+\bar{s}_u u_x-c_1\bar{s}_u-u)>0\\&u\rho \left(B_{u,1}(s)+\bar{s}_u u_x-c_1\bar{s}_u-u\right)+\\&\quad+u\rho^* B_{u,2}(s)+u[\bar{s}_u u_{x,y}-(c_2-\rho c_1)\bar{s}_u-u(a-\rho)]>0\end{split}}.
}
Note that we have the following weak converges in the space $C([-S,T])$ as $u\to\infty$
\bqn{\label{brown_ind}
\begin{aligned}
uB_{u,1}(t)\to B_1(t),\\
uB_{u,2}(t)\to B_2(t),
\end{aligned}
\qquad t\in[-S,T]
}
and further
\bqny{
& &u(\bar{s}_u u_x-c_1\bar{s}_u-u)=u\left[\left(1-\frac{s}{u^2}\right)\left(u+c_1-\frac{x}{u}\right)-c_1\left(1-\frac{s}{u^2}\right)-u\right]\to -s-x,\quad u\to\infty,\\
& &u[\bar{s}_u u_{x,y}-(c_2-\rho c_1)\bar{s}_u-u(a-\rho)]\\
& &=u\left[\left(1-\frac{s}{u^2}\right)\left(u(a-\rho)-\frac{y-\rho x}{u}+c_2-\rho c_1\right)-(c_2-\rho c_1)\left(1-\frac{s}{u^2}\right)-u(a-\rho)\right]\\
& &\to-(a-\rho)s-(y-\rho x),\quad u\to\infty.
}
Consequently
\bqny{
h_u(T,S,x,y)\to h(T,S,x,y),\qquad u\to\infty,
}
where in view of $\eqref{BB}$
\bqny{
h(T,S,x,y)&=&\pk{\exists t\in[0,T]\forall s\in[t-S,t]:\begin{split}&B_1(s)-s-x>0,\\&\rho(B_1(s)-s-x)+\rho^* B_2(s)-(a-\rho)s-(y-\rho x)>0\end{split}}\\
&=&\pk{\exists t\in[0,T]\forall s\in[t-S,t]:\begin{split}&B_1(s)-s>x,\\&\rho B_1(s)+\rho^* B_2(s)-as>y\end{split}}\\
&=&\pk{\exists t\in[0,T]\forall s\in[t-S,t]:\begin{split}&W_1(s)-s>x,\\&W_2(s)-as>y\end{split}}.
}
This convergence is justified by applying continuous mapping theorem for the continuous function $$H_{T,S}(F_1(t),F_2(t))=\sup_{t\in[0,T]}\inf\left(\inf_{s\in[t-S,t]}F_1(t),\inf_{s\in[t-S,t]}F_2(t)\right)$$
and random sequence $(F_{1,x,y,u},F_{2,x,y,u})\in C[-S,T]^2$
\bqny{
F_{1,x,y,u}(s)&=&u(B_{u,1}(s)+\bar{s}_u u_x-c_1\bar{s}_u-u),\\
F_{2,x,y,u}(s)&=&u\rho \left(B_{u,1}(s)+\bar{s}_u u_x-c_1\bar{s}_u-u\right)+u\rho^* B_{u,2}(s)+u[\bar{s}_u u_{x,y}-(c_2-\rho c_1)\bar{s}_u-u(a-\rho)].
}

To finish the proof it is enough to show the dominated convergence as $u\to\infty$ for
\bqny{
I_u(S,T)=\int_{\R^2}h_u(T,S,x,y)\psi_u(x,y)\td x\td y.
}
Note that for $\psi_u(x,y)$ we can write the following upper bound. Fix some $0<\varepsilon<\min(\lambda_1,\lambda_2)$ (such constant exists as in our case both $\lambda_1$ and $\lambda_2$ are greater than zero), and define constants $\lambda_{1,\varepsilon}=\lambda_1+sign(x)\varepsilon$ and $\lambda_{2,\varepsilon}=\lambda_2+sign(y)\varepsilon$. Hence for large enough $u$ and all $x,y\in\R$
\bqn{
\psi_u(x,y)\leq\bar{\psi}:=e^{\lambda_{1,\varepsilon}x+\lambda_{2,\varepsilon}y}.\label{psibar}
}
For $h_u(S,T,x,y)$ we may use Piterbarg inequality (see \cite{Pit96}, Thm 8.1), since 
\bqn{
u^2\E{(B_{u,i}(t)-B_{u,i}(s))^2}<Const|t-s|,\qquad t,s>0\label{PitInequ}
}
for some positive constant and sufficiently large $u$. Thus, for large enough $u$ we have for some positive constant $\bar{C}$
\bqny{
h_u(T,S,x,y)&\leq&\pk{\exists s\in[0,T]:\begin{split}&u(B_{u,1}(s)+\bar{s}_u(u_x-c_1)-u)>0\\&u\rho \left(B_{u,1}(s)+\bar{s}_u( u_x-c_1)-u\right)+u\rho^* B_{u,2}(s)\\&\qquad+u[\bar{s}_u(u_{x,y}-c_2+\rho c_1)-u(a-\rho)]>0\end{split}}\\
&\leq&\bar{h}:=\begin{cases}
\bar{C}e^{-c(x^2+y^2)},~&x,y\geq 0,\\
\bar{C}e^{-cx^2},~&x\geq 0,y<0,\\
\bar{C}e^{-cy^2},~&y\geq 0,x<0,\\
1,~&x,y<0.
\end{cases}
}
Since $\lambda_{1,\varepsilon},\lambda_{2,\varepsilon}$ are positive
\begin{align*}
I_u(S,T)&\leq\int_{\R^2}\bar{h}(T,S,x,y)\bar{\psi}(x,y)\td x\td y & \\
&=\bar{C}\int_{x,y>0}e^{-c(x^2+y^2)+\lambda_{1,\varepsilon}x+\lambda_{2,\varepsilon}y}\td x\td y+\bar{C}\int_{x>0,y<0}e^{-cx^2+\lambda_{1,\varepsilon}x+\lambda_{2,\varepsilon}y}\td x\td y\\
&+\bar{C}\int_{x<0,y>0}e^{-cy^2+\lambda_{1,\varepsilon}x+\lambda_{2,\varepsilon}y}\td x\td y+\bar{C}\int_{x,y<0}e^{\lambda_{1,\varepsilon}x+\lambda_{2,\varepsilon}y}\td x\td y<\infty.
\end{align*}
Hence the proof follows from the dominated convergence theorem.\\
\\
ii) In the case $a\leq\rho$ we use another one transformation
\bqny{
u_x=u+c_1-x/u, \qquad u_y=\rho u +c_2-y.
}
Define also for $u>0$
\bqny{
u_{x,y}:=u_y-\rho u_x=c_2-y-\rho c_1+\rho x/u.
}
In the previous notation
\bqny{
\varphi_{\rho}(u_x,u_y)=:\varphi_{\rho}(u+c_1,\rho u + c_2)\psi_u(x,y),
}
where
\begin{align}
\begin{aligned}
\log \psi_u(x,y)&=(u+c_1,\rho u + c_2)\Sigma^{-1}(x/u,y)^\top-\frac{1}{2}(x/u,y)\Sigma^{-1}(x/u,y)^\top\\
&=(1,\rho)\Sigma^{-1}(x,uy)^\top+(c_1,c_2)\Sigma^{-1}(x/u,y)^\top-\frac{1}{2}(x/u,y)\Sigma^{-1}(x/u,y)\\
&=\frac{1}{1-\rho^2}\left((1,\rho)\begin{pmatrix}1 & -\rho\\-\rho & 1\end{pmatrix}(x,uy)^\top+(c_1,c_2)\begin{pmatrix}1 & -\rho\\-\rho & 1\end{pmatrix}\left(\frac{x}{u},y\right)^\top\right.\\
&\qquad\qquad\qquad\qquad\qquad\qquad\qquad\qquad\qquad\qquad\qquad\left.-\frac{1}{2}\left(\frac{x}{u},y\right)\begin{pmatrix}1 & -\rho\\-\rho & 1\end{pmatrix}\left(\frac{x}{u},y\right)^\top\right)\\
&=\frac{1}{1-\rho^2}\left((1-\rho^2)x+\frac{(c_1-\rho c_2)x}{u}+(c_2-\rho c_1)y-\frac{x^2}{2u^2}-\frac{y^2}{2}+\frac{\rho x y}{u}\right)\\
&\to x-\frac{y^2-2y(c_2-\rho c_1)}{2-2\rho^2}
\end{aligned}\label{iipsi}
\end{align}
as $u\to\infty$. Setting $A(u)=u^{-1}\varphi_\rho(u+c_1,\rho u+c_2)$, we have the following representation for the function $M(u, S, T)$ (write $\bar{s}_u$ for $1-s/u^2$ and recall \eqref{BB})
\bqny{
& &M(u,S,T)=\\
& &\quad=u^{-1}\int_{\R^2}\pk{\exists t\in[\delta(u,T),1]\forall s\in[t,t+S/u^2]:\begin{split}&W_1^*(s)>u\\&W_2^*(s)>au\end{split}\left|\begin{split}&W_1(1)=u_x\\&W_2(1)=u_y\end{split}\right.}\varphi_\rho(u_x,u_y)\td x\td y\\
& &\quad=A(u)\int_{\R^2}\pk{\exists t\in[0,T]\forall s\in[t-S,t]:\begin{split}&W_1^*(\bar{s}_u)>u\\&W_2^*(\bar{s}_u)>au\end{split}\left|\begin{split}&W_1(1)=u_x\\&W_2(1)=u_y\end{split}\right.}\psi_u(x,y)\td x\td y\\
& &\quad=A(u)\int_{\R^2}\pk{\exists t\in[0,T]\forall s\in[t-S,t]:\begin{split}&B_1(\bar{s}_u)-c_1\bar{s}_u>u\\&\rho B_1(\bar{s}_u)+\rho^* B_2(\bar{s}_u)-c_2\bar{s}_u>au\end{split}\left|\begin{split}&B_1(1)=u_x\\&\rho^*B_2(1)=u_{x,y}\end{split}\right.}\psi_u(x,y)\td x\td y\\
& &\quad =:A(u)\int_{R^2}h_u(T,S,x,y)\psi_u(x,y)\td x\td y.
}
Using $B_{u,1}$ and $B_{u,2}$ defined in \eqref{auxBB}
we can represent the function $h_u(T,S,x,y)$ as
\bqny{
h_u(T,S,x,y)=\pk{\exists t\in[0,T]\forall s\in[t-S,t]:\begin{split}&u(B_{u,1}(s)+\bar{s}_u u_x-c_1\bar{s}_u-u)>0\\&u\rho \left(B_{u,1}(s)+\bar{s}_u u_x-c_1\bar{s}_u-u\right)+\\&\quad+u\rho^* B_{u,2}(s)+u[\bar{s}_u u_{x,y}-(c_2-\rho c_1)\bar{s}_u-u(a-\rho)]>0\end{split}}.
}
We have
\bqny{
& &u(\bar{s}_u u_x-c_1\bar{s}_u-u)=u\left[\left(1-\frac{s}{u^2}\right)\left(u+c_1-\frac{x}{u}\right)-c_1\left(1-\frac{s}{u^2}\right)-u\right]\to -s-x,\quad u\to\infty\\
& &u[\bar{s}_u u_{x,y}-(c_2-\rho c_1)\bar{s}_u-u(a-\rho)]=\\
& &\qquad=u\left[\left(1-\frac{s}{u^2}\right)\left(c_2-y-\rho c_1+\frac{\rho x}{u}\right)-(c_2-\rho c_1)\left(1-\frac{s}{u^2}\right)-u(a-\rho)\right]\\
& &\qquad=u\left[\left(1-\frac{s}{u^2}\right)\left(-y+\frac{\rho x}{u}\right)-u(a-\rho)\right]\\
& &\qquad=-u^2(a-\rho)-uy+\rho x+ys/u+\rho x s/u^2.
}
Note that if $a<\rho$, then as $u\to\infty$ the above tends to $+\infty$, and if $a=\rho$ then it tends to $+\infty$ only if $y<0$ and to $-\infty$ if $y>0$. Finally, if $a=\rho$ and $y=0$, then the above tends to $\rho x$.\\
Again using continuous mapping theorem, since \eqref{brown_ind} holds, we have the following convergence (except if $y=0$)
\bqny{
h_u(T,S,x,y)\to h(T,S,x,y),\qquad u\to\infty,
}
where
\bqny{
& &h(T,S,x,y)=\\
& &\quad=\pk{\exists t\in[0,T]\forall s\in[t-S,t]:\begin{split}&B_1(s)-s-x>0,\\&\rho(B_1(s)-s-x)+\rho^* B_2(s)+\infty>0\end{split}}(\mathbb{I}\{a<\rho\}+\mathbb{I}\{a=\rho,~y<0\})\\
& &\quad=\pk{\exists t\in[0,T]\forall s\in[t-S,t]:B_1(s)-s>x}(\mathbb{I}\{a<\rho\}+\mathbb{I}\{a=\rho,~y<0\})\\
& &\quad=\pk{\exists t\in[0,T]\forall s\in[t-S,t]:W_1(s)-s>x}(\mathbb{I}\{a<\rho\}+\mathbb{I}\{a=\rho,~y<0\}).
}
To show the claim we can apply the dominated convergence theorem. Note that for large enough $u$ and all $x,y\in\R$
\bqny{
\log\psi_u(x,y)\leq\bar{\varphi}(x,y)=(1+sign(x)/2)x+\frac{c_2-\rho c_1}{1-\rho^2}y-\frac{y^2}{2}.
}
By Piterbarg inequality (as $\eqref{PitInequ}$ holds here for $i=1$) we can establish that for some positive constant $\bar{C}$
\bqny{
h_u(T,S,x,y)\leq\pk{\exists s\in[0,T]:u(B_{u,1}(s)+\bar{s}_u(u_x-c_1)-u)>0}\leq\bar{h}:=\begin{cases}
\bar{C}e^{-cx^2},~&x\geq 0,\\
1,~&x<0.
\end{cases}
}
Since $(1+sign(x)/2)>0$, then
\bqny{
\int_{\R^2}\bar{h}(x,y)\bar{\varphi}(x,y)\td x \td y<\infty
}
and by the dominated convergence theorem the claim follows.\\
\QED

\begin{remark} In Lemma \ref{L1} if $c_1$, $c_2$ and $S$ depend on $u$, but have finite limits as $u\to\infty$ ($c_1(u)\to c_1^*$, $c_2(u)\to c_2^*$, $S(u)\to S^*$) we can use the same proof. In this case all constants $c_1$ and $c_2$ on the right-hand sides under the density function $\varphi_\rho$ should be left as dependent on $u$, and the rest of them should be replaced by the limited constants $c_1^*$ and $c_2^*$. Also, all constants $S$ on the right-hand sides should be replaced by the limited constant $S^*$. 
\end{remark}

\prooftheo{Th1} Recall first that we define $\delta(u,T)=1-Tu^{-2}$. In view of \nelem{L1} and \eqref{bXh} we
immediately obtain that
\[
\lim_{T\to\infty}\lim_{u\to\infty} \frac{
\pk{ \exists_{t\in [0,\delta(u,T)]}, \forall s\in [t,t+S/u^2]: W_1^*(s)> u, W_2^*(s)> au}
}{P_{0,T}(u,au)} =0.
\]
Hence, using that 
\begin{eqnarray*}
M(u,S,T)\le P_{S/u^2,T}(u,au) \le
\pk{ \exists_{t\in [0,\delta(u,T)]}, \forall s\in [t,t+S/u^2]: W_1^*(s)> u, W_2^*(s)> au} 
+M(u,S,T)
\end{eqnarray*}
we have
\[
\lim_{T\to\infty}\lim_{u\to\infty}
\frac{M(u,S,T)}{P_{S/u^2,T}(u,au)} =1.
\]
Consequently, in view of
Lemma \ref{singA} and the asymptotics of ruin probability of Gaussian vector (see Appendix, Lemma \ref{vec}),
it suffices to prove that
\[
\lim_{T\to\infty}I(S,T)\in (0,\infty),
\]
where $I(S,T)$ is defined in Lemma \ref{singA}. Since $I(S,T) \le I(0,T)$, $I(S,T)$ is growing by $T$ and the finiteness of $\limit{T} I(0,T)$ follows from \cite{mi:18} (see Appendix, Lemma \ref{Ilim}), the claim follows. 
\QED
\\
\\
\auxtheo{Th2} Define $\delta(u,T)=1-Tu^2$ for $T>0$.\\ If $\int_{0}^{1}\mathbb{I}(R_1(t)<0,R_2(t)<0)\td t>L/f(u)$, we have either $$\int_{1-\delta(u,T)}^{1}\mathbb{I}(R_1(t)<0,R_2(t)<0)\td t>L/f(u)$$ or for some point $t_1\in[0,1-\delta(u,T)]$ both $R_1(t_1)$ and $R_2(t_1)$ are lower than zero. On the other hand, if $\int_{1-\delta(u,T)}^{1}\mathbb{I}(R_1(t)<0,R_2(t)<0)\td t>L/f(u)$, then also $\int_{0}^{1}\mathbb{I}(R_1(t)<0,R_2(t)<0)\td t>L/f(u)$. In terms of probabilities it means that
\bqn{
M(u,T)\leq\Psi_L(u,au)\leq M(u,T)+m(u,T),\label{Bconv}
}
where
\bqny{
& &M(u,T)=\pk{\int_{1-\delta(u,T)}^{1}\mathbb{I}(W_1^*(t)>u,W_2^*(t)>au)\td t>L/f(u)},\\
& &m(u,T)=\pk{\exists_{t\in[0,1-\delta(u,T)]}:W_1^*(t)>u,W_2^*(t)>au}.
}
We are going to show further the exact asymptotics of $M(u,T)$ and $m(u,T)=o(M(u,T))$ as $u\to\infty$. It is enough to establish that $\Psi_L(u,au)\sim M(u,T)$ as $u\to\infty$.
Note that for $m(u,T)$, from \cite{mi:18}[Lemma 4.1](see Appendix, Lemma $\ref{inimp}$) we have the following bound for large enough $u$
\bqn{
m(u,T)\leq e^{-T/8}\frac{\pk{W_1^*(1)\ge u, W_2^*(1)\ge au}}{\pk{W_1(1)>\max(c_1,0),W_2(1)>\max(c_2,0)}}.\label{bXh2}
}
\begin{lem}
i) For any $a\in (\rho, 1]$ and any $T>0$ we have as $u\to\infty$
\label{singA}
\BQN
M(u,T)&\sim & u^{-2} \varphi_\rho(u+c_1,au+c_2) I(T),
\EQN
where
$$ I(T)\coloneqq\int_{\R^2}\pk{\int_{0}^{T}\mathbb{I}(W_1(t)-t>x,W_2(t)-at>y)\td t>L}e^{\lambda_1 x+\lambda_2 y}\td x\td y\in(0,\infty).$$
ii) For any $a \le \rho$ and any $T>0$ we have as $u\to\infty$
\BQNY
M(u,T) \sim u^{-1} \varphi_\rho(u+ c_1,\rho u+ c_2) I(T),
\EQNY
where
$$I(T)\coloneqq\int_{\R^2}\pk{\int_0^T\mathbb{I}\left(W_1(t)-t>x\right)\td t>L}\left[\mathbb{I}\{a<\rho\}+\mathbb{I}\{a=\rho,~y<0\}\right]e^{x-\frac{y^2-2y(c_2-c_1\rho)}{2(1-\rho)}}\td x\td y\in(0,\infty).$$
\label{L2.1}
\end{lem}

\prooflem{L2.1}
i) We use the same notation as in Lemma $\ref{L1}$ i). We also have the result $\eqref{ipsi}$.\\
For the function $M(u,T)$ we have
\bqny{
& &M(u,T)\\
& &\quad=u^{-2}\int_{\R^2}\pk{\int_{1-\delta(u,T)}^1\mathbb{I}(W_1^*(t)>u,W_2^*(t)>au)\td t>L/u^2\left|\begin{split}&W_1(1)=u_x\\&W_2(1)=u_y\end{split}\right.}\varphi_\rho(u_x,u_y)\td x\td y\\
& &\quad=A(u)\int_{\R^2}\pk{\frac{1}{u^2}\int_{0}^T\mathbb{I}(W_1^*(\bar{t}_u)>u,W_2^*(\bar{t}_u)>au)\td t>L/u^2\left|\begin{split}&W_1(1)=u_x\\&W_2(1)=u_y\end{split}\right.}\psi_u(x,y)\td x\td y\\
& &\quad=A(u)\int_{\R^2}\pk{\int_0^T\mathbb{I}\left(\begin{split}&B_1(\bar{t}_u)-c_1\bar{t}_u>u\\&\rho B_1(\bar{t}_u)+\rho^* B_2(\bar{t}_u)-c_2\bar{t}_u>au\end{split}\right)\td t>L\left|\begin{split}&B_1(1)=u_x\\&\rho^*B_2(1)=u_{x,y}\end{split}\right.}\psi_u(x,y)\td x\td y\\
& &\quad=:A(u)\int_{R^2}h_u(L,T,x,y)\psi_u(x,y)\td x\td y.
}
Recalling the processes $B_{u,1}$ and $B_{u,2}$ from \eqref{auxBB} we can represent the function $h_u(T,S,x,y)$ as follows
\bqny{
h_u(L,T,x,y)=\pk{\int_0^T\mathbb{I}\left(\begin{split}&u(B_{u,1}(t)+\bar{t}_u u_x-c_1\bar{t}_u-u)>0\\&u\rho \left(B_{u,1}(t)+\bar{t}_u u_x-c_1\bar{t}_u-u\right)+\\&\quad+u\rho^* B_{u,2}(t)+u[\bar{t}_u u_{x,y}-(c_2-\rho c_1)\bar{t}_u-u(a-\rho)]>0\end{split}\right)\td t>L}.
}
Note that we have the same weak converges as in \eqref{brown_ind} and further
\bqn{
\begin{aligned}
&u(\bar{t}_u u_x-c_1\bar{t}_u-u)=u\left[\left(1-\frac{t}{u^2}\right)\left(u+c_1-\frac{x}{u}\right)-c_1\left(1-\frac{t}{u^2}\right)-u\right]\to -t-x,\quad u\to\infty,\\
&u[\bar{t}_u u_{x,y}-(c_2-\rho c_1)\bar{t}_u-u(a-\rho)]\\
&\qquad=u\left[\left(1-\frac{t}{u^2}\right)\left(u(a-\rho)-\frac{y-\rho x}{u}+c_2-\rho c_1\right)-(c_2-\rho c_1)\left(1-\frac{t}{u^2}\right)-u(a-\rho)\right]\\
&\qquad\to-(a-\rho)t-(y-\rho x),\quad u\to\infty.
\end{aligned}\label{constlimit2}
}
Now we want to apply continuous mapping theorem to the function $$H_{T}(F_1,F_2)=\int_0^{T}\mathbb{I}\left(F_1(t)>0,F_2(t)>0\right)\td t$$ and a random sequence $(F_{1,x,y,u},F_{2,x,y,u})\in C([0,T]\to\R^2)$ defined as
\bqny{
F_{1,x,y,u}&=&u(B_{u,1}(t)+\bar{t}_u u_x-c_1\bar{t}_u-u),\\
F_{2,x,y,u}&=&u\rho \left(B_{u,1}(t)+\bar{t}_u u_x-c_1\bar{t}_u-u\right)+u\rho^* B_{u,2}(t)+u[\bar{t}_u u_{x,y}-(c_2-\rho c_1)\bar{t}_u-u(a-\rho)]
}
with exception set $$\Lambda=\{F\in C([0,T]\to\R^2):\mu(F^{-1}(\partial\{(x,y)\in\R^2|x>0,~y>0\}))>0\}.$$ 
First we need to prove that $H_T(F_1,F_2)$ is continuous for $(F_1,F_2)\not\in\Lambda$. Define an area $$\lambda=(F_1,F_2)^{-1}(\partial\{(x,y)\in\R^2|x>0,~y>0\}).$$
For any sequence $(F_{1,n},F_{2,n})$ so that it converges in $C([0,T]\to\R^2)$ to some function $(F_1,F_2)$ as $n\to\infty$ we can define $$(F_{1,n}^\prime(t),F_{2,n}^\prime(t))=\begin{cases}(F_{1,n}(t),F_{2,n}(t)),\quad &t\not\in\lambda,\\
(F_{1}(t),F_{2}(t)),\quad&t\in\lambda.\end{cases}$$
Note that in this case for all $t\in[0,T]$ as $n\to\infty$ $$\mathbb{I}(F_{1,n}^{\prime}>0,F_{2,n}^{\prime}>0)\to\mathbb{I}(F_{1}>0,F_{2}>0).$$
Since $\mu(\lambda)=0$, we have $H_T(F_{1,n}^\prime,F_{2,n}^\prime)=H_T(F_{1,n},F_{2,n})$. Hence, the dominated convergence theorem establishes the continuity of the function $H_T$ at the point $(F_1,F_2)$.\\
From \eqref{brown_ind} and \eqref{constlimit2} we can establish that as $u\to\infty$
\bqny{
F_{1,x,y,u}(t)&\to& B_1(t)-t-x=W_1(t)-t-x,\\
F_{2,x,y,u}(t)&\to& \rho(B_1(t)-t-x)+\rho^* B_2(t)-(a-\rho)t-(y-\rho x)=W_2(t)-at-y.
}
Since $W_1$ and $W_2$ are standard Brownian motions,
\bqny{
\pk{\mu((W_1(\cdot)-\cdot)^{-1}(x))>0}=0,\qquad\qquad \pk{\mu((W_2(\cdot)-a\cdot)^{-1}(y))>0}=0.
}
Consequently, $\pk{(W_1(\cdot)-x-\cdot,W_2(\cdot)-a\cdot-y)\in\Lambda}=0$, and we can apply continuous mapping theorem, which establish that for almost all $L$ positive
\bqny{
h_u(L,T,x,y)\to h(L,T,x,y),\qquad u\to\infty,
}
where
\bqny{
h(L,T,x,y)&=&\pk{\int_0^T\mathbb{I}\left(\begin{split}&B_1(t)-t-x>0\\&\rho(B_1(t)-t-x)+\rho^* B_2(t)-(a-\rho)t-(y-\rho x)>0\end{split}\right)\td t>L}\\
&=&\pk{\int_0^T\mathbb{I}\left(\begin{split}&B_1(t)-t>x\\&\rho B_1(t)+\rho^* B_2(t)-at>y\end{split}\right)\td t>L}\\
&=&\pk{\int_0^T\mathbb{I}(W_1(t)-t>x,W_2(t)-at>y)\td t>L}.
}
To finish the proof it is enough to show the dominated convergence for the integrals
\bqny{
I_u(T)=\int_{\R^2}h_u(L,T,x,y)\psi_u(x,y)\td x\td y.
}
In view of $\eqref{psibar}$ and $\eqref{PitInequ}$ for large enough $u$ we have for some positive constant $\bar{C}$ such that for all $x,y\in\R$
\bqny{
h_u(L,T,x,y)&\leq&\pk{\exists t\in[0,T]:\begin{split}&u(B_{u,1}(t)+\bar{t}_u(u_x-c_1)-u)>0\\&u\rho \left(B_{u,1}(t)+\bar{t}_u( u_x-c_1)-u\right)+u\rho^* B_{u,2}(t)\\&\qquad+u[\bar{t}_u(u_{x,y}-c_2+\rho c_1)-u(a-\rho)]>0\end{split}}\\
&\leq&\bar{h}(T,x,y):=\begin{cases}
\bar{C}e^{-c(x^2+y^2)},~&x,y\geq 0,\\
\bar{C}e^{-cx^2},~&x\geq 0,y<0,\\
\bar{C}e^{-cy^2},~&y\geq 0,x<0,\\
1,~&x,y<0.
\end{cases}
}
Since $\lambda_{1,\varepsilon},\lambda_{2,\varepsilon}$ are positive
\begin{align*}
I_u(T)&\leq\int_{\R^2}\bar{h}(T,x,y)\bar{\psi}(x,y)\td x\td y & \\
&=\bar{C}\int_{x,y>0}e^{-c(x^2+y^2)+\lambda_{1,\varepsilon}x+\lambda_{2,\varepsilon}y}\td x\td y+\bar{C}\int_{x>0,y<0}e^{-cx^2+\lambda_{1,\varepsilon}x+\lambda_{2,\varepsilon}y}\td x\td y\\
&+\bar{C}\int_{x<0,y>0}e^{-cy^2+\lambda_{1,\varepsilon}x+\lambda_{2,\varepsilon}y}\td x\td y+\bar{C}\int_{x,y<0}e^{\lambda_{1,\varepsilon}x+\lambda_{2,\varepsilon}y}\td x\td y<\infty.
\end{align*}
Thus the dominates convergence theorem may be applied and provides us with the claimed assertion. (Note that the constant $I(T)$ is continuous with respect to $L$ (see Appendix, \nelem{MWcont}), so it holds for all $L$ positive)\\
\\
ii) We have the same notation as in Lemma $\ref{L1}$ i).\\
The following representation for the function $M(u,T)$ holds (write $\bar{t}_u$ for $1-t/u^2$ and recall \eqref{BB})
\bqny{
& &M(u,T)=\\
& &\quad=u^{-1}\int_{\R^2}\pk{\int_{1-\delta(u,T)}^T\mathbb{I}(W_1^*(t)>u,W_2^*(t)>au)\td t>L/u^2\left|\begin{split}&W_1(1)=u_x\\&W_2(1)=u_y\end{split}\right.}\varphi_\rho(u_x,u_y)\td x\td y\\
& &\quad=A(u)\int_{\R^2}\pk{\int_0^T\mathbb{I}(W_1^*(\bar{t}_u)>u,W_2^*(\bar{t}_u)>au)\td t>L\left|\begin{split}&W_1(1)=u_x\\&W_2(1)=u_y\end{split}\right.}\psi_u(x,y)\td x\td y\\
& &\quad=A(u)\int_{\R^2}\pk{\int_0^T\mathbb{I}\left(\begin{split}&B_1(\bar{t}_u)-c_1\bar{t}_u>u,\\&\rho B_1(\bar{t}_u)+\rho^* B_2(\bar{t}_u)-c_2\bar{t}_u>au\end{split}\right)\td t>L\left|\begin{split}&B_1(1)=u_x\\&\rho^*B_2(1)=u_{x,y}\end{split}\right.}\psi_u(x,y)\td x\td y\\
& &\quad =:A(u)\int_{R^2}h_u(L,T,x,y)\psi_u(x,y)\td x\td y.
}
Using again $B_{u,1}$ and $B_{u,2}$ as in \eqref{auxBB}
we can represent the function $h_u(L,T,x,y)$ as
\bqny{
h_u(L,T,x,y)=\pk{\int_0^T\mathbb{I}\left(\begin{split}&u(B_{u,1}(t)+\bar{t}_u u_x-c_1\bar{t}_u-u)>0\\&u\rho \left(B_{u,1}(t)+\bar{t}_u u_x-c_1\bar{t}_u-u\right)\\&\quad+u\rho^* B_{u,2}(t)+u[\bar{t}_u u_{x,y}-(c_2-\rho c_1)\bar{t}_u-u(a-\rho)]>0\end{split}\right)\td t>L}.
}
Note that we have the same weak converges as in $\eqref{brown_ind}$. Moreover, we have the convergence $\eqref{iipsi}$.
Using the same arguments as in i) we can use continuous mapping theorem and establish the following convergence for almost all $L$ positive and all $x\in\R,~y\in\R\setminus\{0\}$
\bqny{
h_u(L,T,x,y)\to h(L,T,x,y),\qquad u\to\infty,
}
where
\bqny{
& &h(L,T,x,y)\\
& &\quad=\pk{\int_0^T\mathbb{I}\left(\begin{split}&B_1(t)-t-x>0\\&\rho(B_1(t)-t-x)+\rho^* B_2(t)+\infty>0\end{split}\right)\td t>L}(\mathbb{I}\{a<\rho\}+\mathbb{I}\{a=\rho,~y<0\})\\
& &\quad=\pk{\int_0^T\mathbb{I}\left(B_1(t)-t>x\right)\td t>L}(\mathbb{I}\{a<\rho\}+\mathbb{I}\{a=\rho,~y<0\}).
}
To show the claim we can apply the dominated convergence theorem. Note that for large enough $u$
\bqny{
\log\psi_u(x,y)\leq\bar{\varphi}(x,y)=(1+sign(x)/2)x+\frac{c_2-\rho c_1}{1-\rho^2}y-\frac{y^2}{2}.
}
By Piterbarg inequality we can establish that for some positive constant $\bar{C}$ 
\bqny{
h_u(L,T,x,y)\leq\pk{\exists t\in[0,T]:u(B_{u,1}(t)+\bar{t}_u(u_x-c_1)-u)>0}\leq\bar{h}(x):=\begin{cases}
\bar{C}e^{-cx^2},~&x\geq 0,\\
1,~&x<0.
\end{cases}
}
Since $(1+sign(x)/2)>0$, then
\bqny{
\int_{\R^2}\bar{h}(x)\bar{\varphi}(x,y)\td x \td y<\infty,
}
and by the dominated convergence theorem the claim follows. (As in i), the function $I(T)$ is continuous with respect to $L$, so the assertion holds for all $L\in\R^+$)\\
\QED

\begin{remark} In Lemma \ref{L2.1} if $c_1$, $c_2$ and $L$ depend on $u$, but have finite limits as $u\to\infty$ ($c_1(u)\to c_1^*$, $c_2(u)\to c_2^*$, $L(u)\to L^*$), we can use the same proof. In this case all constants $c_1$ and $c_2$ on the right-hand sides under the density function $\varphi_\rho$ should be left as dependent on $u$, and the rest of them should be replaced by the limited constants $c_1^*$ and $c_2^*$. Also, all constants $L$ on the right-hand sides should be replaced by the limited constant $L^*$. 
\end{remark}

\prooftheo{Th2} Recall that $\delta(u,T)=1-Tu^{-2}$. In view of \nelem{L2.1} and \eqref{bXh2} we
immediately obtain that
\[
\lim_{T\to\infty}\lim_{u\to\infty} \frac{
m(u,T)}
{M(u,T)} =0.
\]
Hence, using \eqref{Bconv}
we have
\[
\lim_{T\to\infty}\lim_{u\to\infty}
\frac{M(u,T)}{\Psi_{L}(u,au)} =1.
\]
Consequently, in view of
Lemma \ref{singA} and the asymptotics of ruin probability of Gaussian vector (see Appendix, Lemma \ref{vec}),
it suffices to prove that
\[
\lim_{T\to\infty}I(T)\in (0,\infty),
\]
where $I(T)$ is defined in Lemma \ref{singA}. Since $I(T) \le I(L,T)$ defined in Lemma \ref{L1}, $I(T)$ is growing by $T$ and $\limit{T} I(L,T) < \IF$, the claim follows. 
\QED

\proofprop{Th3}
Using the formula of conditional probability and recall the self-similarity of Brownian motion
\bqny{
\pk{u^2(1-\tau_{L_1}(u))\geq x|\tau_{L_2}(u)\leq 1}&=&\frac{\pk{\tau_{L_1}(u)\leq 1-x/u^2}}{\pk{\tau_{L_2}(u)\leq 1}}\\
&=&\frac{\pk{{\displaystyle \int_0^{1-x/u^2}}\mathbb{I}\left(\begin{array}{ccc}W_1^*(t)>u\\W_2^*(t)>au\end{array}\right)\td t>L_1/u^2}}{\pk{{\displaystyle \int_0^{1}}\mathbb{I}\left(\begin{array}{ccc}W_1^*(t)>u\\W_2^*(t)>au\end{array}\right)\td t>L_2/u^2}}\\
&=&\frac{\pk{{\displaystyle \int_0^{1}}\mathbb{I}\left(\begin{array}{ccc}W_1^*((1-x/u^2)t)>u\\W_2^*((1-x/u^2)t)>au\end{array}\right)\td (1-x/u^2)t>L_1/u^2}}{\pk{{\displaystyle \int_0^{1}}\mathbb{I}\left(\begin{array}{ccc}W_1^*(t)>u\\W_2^*(t)>au\end{array}\right)\td t>L_2/u^2}}\\
&=&\frac{\pk{{\displaystyle \int_0^{1}}\mathbb{I}\left(\begin{array}{ccc}W_1(t)>\frac{u}{\sqrt{1-x/u^2}}+c_1\sqrt{1-x/u^2}t\\W_2(t)>a\frac{u}{\sqrt{1-x/u^2}}+c_2\sqrt{1-x/u^2}t\end{array}\right)\td t>\frac{L_1/(1-x/u^2)^2}{(u/\sqrt{1-x/u^2})^2}}}{\pk{{\displaystyle \int_0^{1}}\mathbb{I}\left(\begin{array}{ccc}W_1^*(t)>u\\W_2^*(t)>au\end{array}\right)\td t>L_2/u^2}}.
}
Using Theorem $\ref{Th2}$
\bqny{
\frac{\pk{{\displaystyle \int_0^{1}}\mathbb{I}\left(\begin{array}{ccc}W_1(t)>\frac{u}{\sqrt{1-x/u^2}}+c_1\sqrt{1-x/u^2}t\\W_2(t)>a\frac{u}{\sqrt{1-x/u^2}}+c_2\sqrt{1-x/u^2}t\end{array}\right)\td t>\frac{L(1-x/u^2)^2}{(u/\sqrt{1-x/u^2})^2}}}{\pk{{\displaystyle \int_0^{1}}\mathbb{I}\left(\begin{array}{ccc}W_1^*(t)>u\\W_2^*(t)>au\end{array}\right)\td t>L/u^2}}\qquad\qquad\qquad\\
\sim\Gamma(L_1,L_2)\frac{\pk{W_1(1)>\frac{u}{\sqrt{1-x/u^2}}+c_1\sqrt{1-x/u^2},W_2(1)>a\frac{u}{\sqrt{1-x/u^2}}+c_2\sqrt{1-x/u^2}}}{\pk{W_1^*(1)>u,W_2^*(1)>au}}.
}
Notice that (write $\varphi_\rho(x,y)$ for the pdf of vector $(W_1(1),W_2(1))$)
\bqny{
\varphi\left(\frac{u}{\sqrt{1-x/u^2}}+c_1\sqrt{1-x/u^2},\frac{u}{\sqrt{1-x/u^2}}a+c_2\sqrt{1-x/u^2}\right)=\varphi(u+c_1,au+c_2)\psi_u^*(a,c_1,c_2),
}
where
\bqny{
\lim_{n\to\infty}\log\psi_u^*(\rho,c_1,c_2)=-x\frac{1-2a\rho+a^2}{2-2\rho^2},
}
hence by Lemma \ref{vec} the claim follows if $a>\rho$ . For the case $a\leq\rho$ notice that 
\bqny{
\varphi\left(\frac{u}{\sqrt{1-x/u^2}}+c_1\sqrt{1-x/u^2},\frac{u}{\sqrt{1-x/u^2}}\rho+c_2\sqrt{1-x/u^2}\right)=\varphi(u+c_1,\rho u+c_2)\psi_u^*(\rho,c_1,c_2),
}
where
\bqny{
\lim_{u\to\infty}\log\psi^*_u(\rho,c_1,c_2)=-x/2.
}
This finishes the proof in the case $a\leq\rho$ again using Lemma \ref{vec}
\QED

\section{Appendix}
\subsection{Parisian ruin for non-collapsing interval}
Consider now a general probability $p_{H,T}(u,au)$ with some fixed constant $H$. As it was mentioned in $(\ref{4})$, we can reduce the probability by putting $t=T$ instead of supremum.
\bqny{ \pk{ \forall t\in [T,T+H]: W_1^*(t)>u, W_2^*(t) > au} \le p_{H,T}(u,au).}
We can present $W_2(t)$ using the correlation coefficient $\rho$ as $\rho W_1(t)+\rho^* B(t)$, where $\rho^*=\sqrt{1-\rho^2},$ and $B(t)$ is an independent Brownian motion. Note that if $W_1^*(t)>u$ and $B(t)>(a-\rho)u+(c_2-\rho c_1)t$,
then also $W_2^*(t)>au$. Since $W_1$ and $B$ are independent
\bqny{
\pk{\forall t\in[T,T+H]:R_1(t)<0}\pk{\forall t\in[T,T+H]:B(t)>(a-\rho)u+(c_2-\rho c_1)t}\leq p_{H,T}(u,au).
}
In case $\rho >0$ and $\rho>a$, the probability $\pk{\forall t\in[T,T+H]:B(t)>(a-\rho)u+(c_2-\rho c_1)t}$ tends to one when $u$ tends to infinity. So, for any positive $\varepsilon$ for large enough $u$ we derived the following lower bound
\bqny{
(1-\varepsilon)\pk{\forall t\in[T,T+H]: W_1(t)-c_1t>u}\leq p_{H,T}(u,au). 
}
To find an upper bound we can simply put $H=0$ and forget about the inequality for $W_2$:
\bqny{
p_{H,T}(u,au)\leq\pk{\sup_{t \in[0,T]}(W_1(t)-c_1 t)>u}.
}
\\
As it was mentioned, the received upper bound is asymptotically equal to $\pk{W_1^*(T)>u}$, and the lower bound is asymptotically equal to $\pk{W_1^*(T)>u}/u$ as u tends to infinity. 
\subsection{Auxiliary lemmas} The first Lemma presents an asymptotics of ruin probability of Gaussian vector
\begin{lem} \label{vec}Let $X_1$ and $X_2$ be Gaussian random values with correlation $\rho\in(-1,1)$. Let also $c_1,c_2$ be two given real constants and $a \le 1$ be given. Write further $\varphi_{\rho}(x,y)$ for the joint density function of vector $(X_1,X_2)$.\\
i) If $a\in (\rho,1]$, then as $u\to \IF$
\bqny{
\pk{X_1>u+c_1,X_2>au+c_2} \sim \frac{u^{-2}}{\lambda_1\lambda_2}\varphi_\rho(u+c_1, a u+ c_2),
}
where
\bqny{
\lambda_1= \frac{1- a\rho}{1- \rho^2}, \quad \lambda_2=\frac{a- \rho }{1- \rho^2}.
}
ii) If $a\le \rho$, then we have as $u\to \IF$
\bqny{
\pk{X_1>u+c_1,X_2>au+c_2} \sim \sqrt{2 \pi (1- \rho^2)} \Phi^*(c_1 \rho- c_2) e^{ \frac{ (c_2- \rho c_1)^2}{2(1- \rho^2)}} u^{-1}\varphi_\rho(u+c_1, \rho u+ c_2) ,
}
where $\Phi^*(c_1 \rho- c_2)=1$ if $a< \rho$ and $\Phi^*$ is the df of $\sqrt{1- \rho^2} X_1
$ when $a=\rho$.\\
\end{lem}

\prooflem{vec}
i) Using dominated convergence theorem, as $u\to\infty$
\bqny{
\pk{X_1>u+c_1,X_2>au+c_2}&=&\int_{x,y>0}\varphi_{\rho}(u+c_1+x,au+c_2+y)\td x\td y\\
&=&\varphi_\rho(u+c_1,au+c_2)\int_{x,y>0}e^{-u\left(\lambda_1x+\lambda_2 y\right)}\frac{\varphi_\rho(x+c_1,y+c_2)}{\varphi_\rho(c_1,c_2)}\td x\td y\\
&=&\frac{\varphi_\rho(u+c_1,au+c_2)}{u^2}\int_{x,y>0}e^{-\lambda_1x-\lambda_2y}\frac{\varphi_\rho(c_1+x/u,c_2+y/u)}{\varphi_\rho(c_1,c_2)}\td x\td y\\
&\sim&\frac{\varphi_\rho(u+c_1,au+c_2)}{u^2}\int_{x,y>0}e^{-\lambda_1x-\lambda_2y}\td x\td y\\
&=&\frac{\varphi_\rho(u+c_1,au+c_2)}{\lambda_1\lambda_2u^2}.
}
ii)Again using dominated convergence theorem as $u\to\infty$ (denote $C=0$ if $a=\rho$ and $C=-\infty$ otherwise)
\bqny{
\pk{X_1>u+c_1,X_2>au+c_2}&=&\int\limits_{\substack{x>0\\ y>(a-\rho)u}}\varphi_{\rho}(u+c_1+x,\rho u+c_2+y)\td x\td y\\
&=&\varphi_{\rho}(u+c_1,\rho u+c_2)\int\limits_{\substack{x>0\\ y>(a-\rho)u}}e^{-ux}\frac{\varphi_{\rho}(c_1+x,c_2+y)}{\varphi_{\rho}(c_1,c_2)}\td x\td y\\
&=&\frac{\varphi_{\rho}(u+c_1,\rho u+c_2)}{u}\int\limits_{\substack{x>0\\ y>(a-\rho)u}}e^{-x}\frac{\varphi_{\rho}(c_1+x/u,c_2+y)}{\varphi_{\rho}(c_1,c_2)}\td x\td y\\
&\sim&\frac{\varphi_{\rho}(u+c_1,\rho u+c_2)}{u}\int\limits_{\substack{x>0\\ y>C}}e^{-x}\frac{\varphi_{\rho}(c_1,c_2+y)}{\varphi_{\rho}(c_1,c_2)}\td x\td y\\
&=&\frac{\varphi_{\rho}(u+c_1,\rho u+c_2)}{u}\int\limits_{y>C}\frac{\varphi_{\rho}(c_1,c_2+y)}{\varphi_{\rho}(c_1,c_2)}\td y\\
&=&\frac{\varphi_{\rho}(u+c_1,\rho u+c_2)}{u}\int\limits_{y>C}e^{-\frac{1}{2}\left(\frac{y^2}{1-\rho^2}+2\frac{c_2-\rho c_1}{1-\rho^2}y\right)}\td y\\
&=&\frac{\varphi_{\rho}(u+c_1,\rho u+c_2)}{u}e^{\frac{(c_2-\rho c_1)^2}{1-\rho^2}}\int\limits_{y>C}e^{-\frac{1}{2}\left(\frac{(y+c_2-\rho c_1)^2}{1-\rho^2}\right)}\td y\\
&=&\frac{\varphi_{\rho}(u+c_1,\rho u+c_2)}{u}e^{\frac{(c_2-\rho c_1)^2}{1-\rho^2}}\sqrt{2\pi(1-\rho^2)}\int\limits_{y<\rho c_1-c_2-C}\frac{e^{-\frac{1}{2}\frac{y^2}{1-\rho^2}}}{\sqrt{2\pi(1-\rho^2)}}\td y\\
&=&\frac{\varphi_{\rho}(u+c_1,\rho u+c_2)}{u}e^{\frac{(c_2-\rho c_1)^2}{1-\rho^2}}\sqrt{2\pi(1-\rho^2)}\Phi^*(\rho c_1-c_2).
}
\QED

\begin{remark} Lemma \ref{vec} may be used also if $c_1$ and $c_2$ depend on $u$, but have finite limits as $u\to\infty$ ($c_1(u)\to c_1^*$, $c_2(u)\to c_2^*$). In this case all constants $c_1$ and $c_2$ under the density function $\varphi_\rho$ and under the probability signs should be left as dependent on $u$, and the rest of them should be replaced by the limited constants $c_1^*$ and $c_2^*$.
\end{remark}
\begin{lem}\label{MWcont}
Let $X_1(t),X_2(t)$ for $t\geq0$ have representation \eqref{BB}. Let also $\lambda_1,\lambda_2,a,T$ be positive constants and $c_1,c_2$ be real constants. Then the functions
\bqny{
I_1(L)&=&\int_{\R^2}\pk{\int_{0}^{T}\mathbb{I}(X_1(t)-t>x,X_2(t)-at>y)\td t>L}e^{\lambda_1 x+\lambda_2 y}\td x\td y,\\
I_2(L)&=&\int_{\R^2}\pk{\int_0^T\mathbb{I}\left(X_1(t)-t>x\right)\td t>L}\left[\mathbb{I}\{a<\rho\}+\mathbb{I}\{a=\rho,~y<0\}\right]e^{x-\frac{y^2-2y(c_2-c_1\rho)}{2(1-\rho)}}\td x\td y
}
are continuous for $L\in(0,+\infty)$.
\end{lem}
\prooflem{MWcont}
Consider the function $I_1(L)$. The proof for $I_2(L)$ will be the same. To show the continuity of $I_1(L)$ it is enough to show
\bqny{
I_1^*(L):=\int_{\R^2}\pk{\int_{0}^{T}\mathbb{I}(W_1(t)-t>x,W_2(t)-at>y)\td t=L}e^{\lambda_1 x+\lambda_2 y}\td x\td y=0
}
for all positive $L$. Fix some $L>0$ and let
\bqny{
A_{x,y}=\left\{f_1,f_2\in C[0,T]:\int_{0}^{T}\mathbb{I}(f_1(t)-t>x,f_2(t)-at>y)\td t=L\right\}.
}
First note that for any fixed $y_0\in\R$ the sets $A_{x_1,y_0}$ and $A_{x_2,y_0}$ are non-overlapping for $x_1\not= x_2$. Define
$$\mathcal{X}=\{x\in\R: \pk{A_{x,y_0}}>0\} $$
and sets
$$\mathcal{X}_n=\{x\in\R: \pk{A_{x,y_0}}>1/n\} $$
Since $A_{x,y_0}$ are non-overlapping for different $x\in\R$, $|\mathcal{X}_n|<n$. In addition, $\mathcal{X}=\cup_{n\in\N}\mathcal{X}_n$.
Thus, the set $\mathcal{X}$ is countable, establishing the proof.

\QED

The following lemmas were already established in \cite{mi:18}:
\begin{lem}\label{Ilim} Let $W_1(t)$, $W_2(t)$ for $t\geq 0$ have representation \eqref{BB} Assume also some real constants $c_1$, $c_2$, and a real constant $a\leq 1$.\\
i) If $a\in(\rho,1]$, then for any positive $T$, $\lambda_1$, $\lambda_2$
\bqny{
I(T)=\int_{\R^2}\pk{\sup_{t\in[0,T]}\min\left(W_1(t)-t-x,W_2(t)-at-x\right)>0}e^{\lambda_1 x+\lambda_2 y}\td x\td y
}
is positive and bounded by some positive constant $M$ which does not depend on $T$.\\
ii) If $a\leq\rho$, then for any positive $T$
\bqny{
I(T)=\int_{\R^2}\pk{\sup_{t\in[0,T]}\left(W_1(t)-t-x\right)>0}\biggl[\mathbb{I}(a<\rho) + \mathbb{I}(a=\rho, y<0)\biggr]e^{x-\frac{y^2-2y(c_2-c_1\rho)}{2(1-\rho^2)}}\td x\td y
}
is positive and bounded by some positive constant $M$ which does not depend on $T$.
\end{lem}

\begin{lem}\label{inimp}
For any $T>0$, $a\in(-\infty,1]$ and sufficiently large u
\bqny{
\pk{\exists_{t\in[0,\delta(u,T)]}:W_1^*(t)>u, W_2^*(t)>au}\leq e^{-T/8}\frac{\pk{W_1^*(1)\geq u,W_2^*(1)\geq au}}{\pk{W_1(1)>\max(c_1,0),W_2(1)>\max(c_2,0)}},
}
where $W_1(t)$ and $W_2(t)$ for $t\geq 0$ are two standard Brownian motions satisfy \eqref{BB} and $\delta(u,T)=1-T/u^2$.
\end{lem}

\bibliographystyle{ieeetr}
\def\polhk#1{\setbox0=\hbox{#1}{\ooalign{\hidewidth
  \lower1.5ex\hbox{`}\hidewidth\crcr\unhbox0}}}

\end{document}